\documentclass[13pt]{article}
\include{defs}
\usepackage{latexsym}
\usepackage{amsmath,enumerate}
\usepackage{pstricks}
\usepackage{amssymb}
\usepackage{latexsym}
\usepackage{amscd}
\usepackage{amsmath}
\usepackage{eufrak}
\usepackage{txfonts}
\usepackage{indentfirst}
\newtheorem{thm}{Theorem}
\newtheorem{lem}{Lemma}

\newtheorem{defn-lem}[thm]{Definition-Lemma}

\newtheorem{prop}[thm]{Proposition}

\def\wt{\widetilde}

\def \A{{\mathbb A}}
\allowdisplaybreaks[4]

\def \A{{\mathbb A}}
\def \C{{\mathbb C}}

\def \G{{\bf G}}

\def \GSp{{\mathrm {GSp}}}

\def \GSO{{\mathrm {GSO}}}
\def \GO{{\mathrm {GO}}}
\def \GL{{\mathrm {GL}}}
\def \PGL{{\mathrm {PGL}}}
\def \PGO{{\mathrm {PGO}}}
\def \PGSp{{\mathrm {PGSp}}}
\def \Sp{{\mathrm {Sp}}}
\def \SL{{\mathrm {SL}}}
\def\ds{ \displaystyle}

\def\az{\alpha}

\def\oz{\omega}

\def\sz{\sigma}

\def\lz{\lambda}

\begin{document}

\begin{center}

{\Large \bf{Endoscopy and the Transfer from $\GSp (4)$ to $\GL (4)$}}

\end{center}

\bigskip

\begin{center}

{\Large {Bogume Jang}}

\end{center}

\bigskip

\begin{enumerate}

\item[{\large \bf{1}}  ] {\large\bf{Introduction}}

\bigskip

\noindent The purpose of this work is to prove Jiang's conjecture\cite{J} based on the analysis under the assumption of the functorial liftings and the endoscopic liftings so that we can see the existence of a L-function of a cuspidal representation of $ \GSp (4,\A ) \times \GSp (4,\A )$ which has a pole of order 2 at $s=1$, even for globally generic representations.

In \cite{J}, Dihua Jiang studies the degree 16 Rankin product L-function for $\GSp (4) \times \GSp (4)$, where $\GSp (4)$ is the reductive group of symplectic similitudes of rank 2. More precisely, this L-function is defined as fo1lows: Let $\pi _1$ and $\pi _2$ be irreducible automorphic cuspidal representations of $\GSp (4,\A)$ with trivial central characters and let $\rho $ be the standard representation of $\GSp (4,\C )$, the complex dual group of $\GSp (4)$ \cite{B76}. The degree 16 standard L-function is $L^S(s, \pi _1 \otimes \pi _2,\rho \otimes \rho )$.

As explained by Jiang\cite{J} the following commutative diagram with L-homomorphisms: $\sigma _2=\sigma \circ \sigma _1$ and

$\begin{array}[c]{ccccc}
\Sp (4,\C) &             & \stackrel{\sigma }{\rightarrow} & & \SL (4,\C)\\
         & \nwarrow ^{\sigma _1} & &  & {}^{\sigma _2}\nearrow \\
 & & \SL (2,\C) \times \SL (2,\C) & &
\end{array}$

will allow us to predict analytic properties of the L-function $L^S(s, \pi _1 \otimes \pi _2,\rho \otimes \rho )$ for $Re(s)>0$ as follows:

(1) If neither $\pi _1$ nor $\pi _2$ is an endoscopic lifting via $\sigma _1$, then
\[L^S(s,\pi_1 \otimes \pi_2,(\rho  \otimes \rho )\circ \sigma )=L^S(s,\sigma (\pi_1) \otimes \sigma (\pi_2), \rho  \otimes \rho )\]
is holomorphic for all $s$ except at $s=1$ where the L-function $L^S(s,\pi _1 \otimes \pi _2,\rho  \otimes \rho )$ has a simple pole if and only if $\sigma (\pi_2) = \sigma (\pi_1)^\vee$, the contragredient representation of $\sigma (\pi_1)$.

(2) If only one of $\pi _1$ and $\pi _2$ is an endoscopic lifting via $\sigma _1$, then the L-function $L^S(s,\pi _1 \otimes \pi _2,\rho  \otimes \rho )$ is holomorphic for all $s$.

In fact, if, say $\pi _1=\sigma _1 (\pi _1 ^{(1)} \otimes \pi _2^{(2)})$, an endoscopic lifting via $\sigma _1$, then one has
\[\sigma (\pi _1)= \sigma _2(\pi _1 ^{(1)} \otimes \pi _2^{(1)})=\pi _1 ^{(1)} \oplus \pi _2^{(1)}~(\mathrm{automorphic~induction}).  \]
Thus the L-function $L^S(s,\pi _1 \otimes \pi _2,(\rho  \otimes \rho )\circ \sigma )$ has following properties:

$L^S(s,\pi _1 \otimes \pi _2,(\rho  \otimes \rho )\circ \sigma )$

=$L^S(s,\sigma (\pi _1) \otimes \sigma (\pi _2),\rho  \otimes \rho )$

=$L^S(s,(\pi _1^{(1)} \oplus \pi _2^{(1)})\otimes \sigma (\pi _2),\rho  \otimes \rho )$

=$L^S(s,\pi _1^{(1)} \otimes \sigma (\pi _2),\rho  \otimes (\rho  \circ \sigma _2))\cdot L^S(s,\pi _2^{(1)} \otimes \sigma (\pi _2),\rho  \otimes (\rho  \circ \sigma _2)).$

Since $L^S(s,\pi _1^{(1)} \otimes \sigma (\pi _2),\rho  \otimes (\rho  \circ \sigma _2))$ and $L^S(s,\pi _2^{(1)} \otimes \sigma (\pi _2),\rho  \otimes (\rho  \circ \sigma _2))$ are L-functions of $\PGL (2) \times \PGL (4)$, they are holomorphic for all $s$. Thus the product of these two L-functions is holomorphic for all $s$.

(3) If both of $\pi _1$ and $\pi _2$ are endoscopic liftings by means of $\sigma _1$, then the L-function $L^S(s,\pi _1 \otimes \pi _2,\rho  \otimes \rho )$ is holomorphic for all $s$ except for $s=1$ where the L-function may achieve a pole of degree at most two, according to the following discussion:

Assume that $\pi _1=\sigma _1(\pi _1 ^{(1)} \otimes \pi _2^{(1)})$ and $\pi _2=\sigma _1(\pi _1 ^{(2)} \otimes \pi _2^{(2)})$ are endoscopic liftings via $\sigma _1$, then one has

$L^S(s,\pi _1 \otimes \pi _2,(\rho  \otimes \rho )\circ \sigma )$

=$L^S(s,\sigma (\pi _1) \otimes \sigma (\pi _2))$

=$L^S(s,(\pi _1^{(1)} \oplus \pi _2^{(1)})\otimes (\pi _1^{(2)} \oplus \pi _2^{(2)}))$

=$L^S(s,\pi _1^{(1)} \otimes \pi _1^{(2)})\cdot L^S(s,\pi _1^{(1)} \otimes \pi _2^{(2)})\cdot L^S(s,\pi _2^{(1)} \otimes \pi _1^{(2)})\cdot L^S(s,\pi _2^{(1)} \otimes \pi _2^{(2)}).$

Each of the L-function $L^S(s,\pi _i^{(1)} \otimes \pi _j^{(2)})$ is a standard L-function of $\PGL (2) \times \PGL (2)$, which is holomorphic except for $s=1$ where the L-funtion $L^S(s,\pi _i^{(1)} \otimes \pi _j^{(2)})$ has a simple pole if and only if $\pi _i^{(1)}$ is the contragredient representation of $\pi _j^{(2)}$. Note that automorphic cuspidal representations of $\PGL (2)$ are self dual. Thus one has following cases:

(3a) If $L^S(s,\pi _1 \otimes \pi _2,\rho  \otimes \rho )$ has a pole at $s=1$ of fourth degree, then each of $L^S(s,\pi _i^{(1)} \otimes \pi _j^{(2)})$ has a simple pole at $s=1$. Thus all these four representations are equivalent to each other. This implies that $\pi _1=\pi _2=\sigma _1(\pi \otimes \pi )$ for a cusp form $\pi $ of $\PGL (2)$. According to Rallis' theory of tower of theta liftings, the first occurrence of the theta lifting of the automorphic cuspidal representation $\pi \otimes \pi $ of $\PGO (2,2)(\PGL (2) \times \PGL (2))$ is on the group $\PGSp  (2)$. Therefore $\sigma _1(\pi \otimes \pi )$ is no longer a cusp form on $\PGSp (4)$. It follows that the degree 16 standard L-function $L^S(s,\pi _1 \otimes \pi _2,\rho  \otimes \rho )$ of $\GSp (4) \times \GSp (4)$ can not have a pole at $s=1$ of degree greater than three. Note that the L-function $L^S(s,\pi _1 \otimes \pi _2,\rho  \otimes \rho )$ can not have a pole at $s=1$ of degree three following the same argument.

(3b) If $L^S(s,\pi _1 \otimes \pi _2,\rho  \otimes \rho )$ has a pole at $s=1$ of degree two, the only case that both $\pi _1=\sigma _1(\pi _1 ^{(1)} \otimes \pi _2^{(1)})$ and $\pi _2=\sigma _1(\pi _1 ^{(2)} \otimes \pi _2^{(2)})$ are cusp forms on $\PGSp (4)$ is $\pi _1 ^{(1)} = \pi _1^{(2)}$ and $\pi _2 ^{(1)} = \pi _2^{(2)}$ with $\pi _1 ^{(1)} \neq \pi _2^{(1)}$. This implies that $\pi _2$ is the contragredient representation of $\pi _1$ and $\pi _1=\sigma _1(\pi _1 ^{(1)} \otimes \pi _2^{(1)})$ is a cusp form in the image of theta lifting from $\PGO (2,2)$ to $\PGSp (4)$.

(3c) If $L^S(s,\pi _1 \otimes \pi _2,\rho  \otimes \rho )$ has only a simple pole at $s=1$, then one has $\pi _1 ^{(1)} \neq \pi _2^{(1)} \neq \pi _2 ^{(2)} \neq \pi _1^{(1)}$ and $\pi _1 ^{(1)} = \pi _1^{(2)}$. It follows that $\pi _1=\sigma _1(\pi _1 ^{(1)} \otimes \pi _2^{(1)})$ and $\pi _2=\sigma _1(\pi _1 ^{(1)} \otimes \pi _2^{(2)})$ are cusp forms on $\PGSp (4)$, which are not contragredient to each other.

Based on the above analysis under the assumption of the functorial liftings and the endoscopic liftings, Jiang\cite{J} predicts the following conclusions:

\noindent {\bf  { Conjecture }}

(*) For a generic cusp form $\pi $ on $\GSp (4, \A)$ with trivial central character, the automorphic L-function $L^S(s, \pi \otimes \pi ^{\vee}, \rho  \otimes \rho )$ is holomorphic for all $s$ except for $s=1$ where the L-function $L^S(s, \pi \otimes \pi ^{\vee}, \rho  \otimes \rho )$ has a pole of degree at most two.

(**) For a generic cusp form $\pi $ on $\GSp (4,\A )$ with trivial central character, the automorphc L-function $L^S(s, \pi \otimes \pi ^{\vee}, \rho  \otimes \rho )$ achieves the second degree pole at $s=1$ if and only if the generic cusp form $\pi $ is a nonzero endoscopic lifting of a generic cusp form on $\PGL (2, \A ) \times \PGL (2, \A )$.

Part (3b) shows the existence of a L-function of a cuspidal representation of $ \GSp (4,\A) \times \GSp (4,\A)~$which has a pole of order 2 at $s=1$, even for globally generic representations. Recent work of Asgari and Shahidi \cite{AS} has made it possible to prove these statements and get the full analytic continuation of this L-function and more. This occupies the bulk of this paper which we now explain.

\bigskip

\item[{\large\bf{2}}  ] {\large\bf{Classification Theory}}

\bigskip

\noindent In this section, we will see the classification theorem to show that the transferred representation from GSp(4) to GL(4) is the isobaric sum of representations in GL(2)'s.

Let $F$ be a global field and $\A$ be the ring of adeles. 	

We can obtain a classification theorem for automorphic forms on $\mathrm{GL}(r)$ which is a precise analogue for this group of the known results for
	local groups by \cite{JS}.

	Accordingly let $\mathrm{P}~$be a standard parabolic subgroup of $\mathrm{GL}(r)$ of type $(r_1,~r_2, \cdots,~r_u)$. The quotient of $\mathrm{P}$ with its unipotent radical
	$\mathrm{U}_P~$is isomorphic to the group
	\[ \mathrm{M=\GL}(r_1) \times \mathrm{\GL} (r_2) \times \cdots \times \mathrm{\GL} (r_u).\]

	For each $j$, $1 \leq j \leq u$, let $\sz_j~$be an automorphic cuspidal representation of $\mathrm{GL} (r_j, \mathbb{A})$. For each place $\nu~$the
	representation $\sigma _{\nu}=\otimes _j \sigma_{j \nu}~$of the group $\mathrm{M}(F_{\nu})~$can be regarded as a representation of $\mathrm{P}(F_{\nu})$ trivial on $\mathrm{U}(F_{\nu})$, it induces an admissible representation of $\mathrm{GL}(r,F_{\nu})~$
	which we will denote by
	\[ \xi _{\nu}= \mathrm{Ind}(\mathrm{GL}(r,F_{\nu}),~\mathrm{P}(F_{\nu});~\sigma_{\nu}).\]

	One obtains a family of irreducible admissible representations of $\mathrm{GL}(r,\mathbb{A})~$by taking for each irreducible component
	$\pi _{\nu}~$of the representation $\xi _{\nu}~$and forming the tensor product $\pi= \otimes _{\nu} \pi_{\nu}$. On the other hand, with $ \sigma =
	\sigma _1 \otimes \sigma _2 \otimes \cdots \otimes \sigma _u$, one can define globally an induced representation
	\[ \xi = \mathrm{Ind} (\mathrm{GL}(r,\mathbb{A}),~\mathrm{P}(\mathbb{A});~\sigma).\]

	Of course $ \xi = \otimes _{\nu} \xi _{\nu}$.

	Let $\mathrm{Q}~$be another standard parabolic say of type $(s_1,~s_2, \cdots ,~s_w)~$and $ \tau _j~$an automorphic cuspidal representation of
	$\mathrm{GL} (s_j, \mathbb{A})$. As before let $\tau _{\nu}= \otimes _j \tau _{j_{\nu}}$,
	\[ \eta _{\nu} = \mathrm{Ind} (\mathrm{GL}(r,F_{\nu}),~\mathrm{Q}(F_{\nu}); \tau _{\nu}),\]
	and
	\[ \eta = \mathrm{Ind} (\mathrm{GL}(r,\mathbb{A}),~\mathrm{Q}(\mathbb{A}); \tau),\]
	where $ \tau = \tau _1 \otimes \tau _2 \otimes \cdots \otimes \tau _w$. We may ask whether $ \xi ~$and $ \eta ~$have a common constituent. Suppose
	$\mathrm{P}~$and $\mathrm{Q}~$are associate and there is a permutation $ \phi~$of $\{1,~2,~ \cdots ,~u \}~$such that $s_j=r_{ \phi (j)}$. Suppose
	moreover that $ \tau _j = \sigma_{\phi (j)}$. We will say in this situation that the pairs $( \sigma ,~\mathrm{P})~$and $( \tau ,~\mathrm{Q})~$are
	associate.
	When this is so the representations $ \xi _{\nu}~$and $ \eta _{\nu} ~$have the same character, and therefore the same components. In particular if both $ \xi _{\nu}~$and $ \eta _{\nu} ~$are unramified then their unique unramified components are the same. In other
	words the irreducible components of $ \xi ~$and $ \eta~$are the same.

    And the converse is also true by theorem 4.4 in \cite{JS}.
    \begin{prop}
    Let P, Q, $\sz _j$ and $\tau _k$ be as above. Let $S$ be a finite set of places containing all the places at infinity. Suppose that for $\nu \notin S$ the representations $\sz _{j \nu }$ and $\tau _{k \nu }$ are unramified and that the representations $\xi _{\nu }$ and $\eta _{\nu }$ of $\mathrm{GL}(r,F _{\nu })$ they induce have the same unramified component. Then the pairs $(\sz , P)$ and $(\tau , Q )$ are associate.
    \end{prop}

\bigskip

\item[{\large\bf{3}}  ] {\large\bf{Transfer from GSO(4) to GL(4)}}

\bigskip


In this section, we will see the relation between $\GSO (4)$ and $\GL (4)$ which we will need later to show the existence of a representation transferred from GSO(4) to GL(4).

Let $k$ be a number field with algebraic closure $\bar{k}$. Let $V$ be a finite dimensional vector space over $k$ equipped with a non-degenerate symmetric bilinear form $B:V \times V \rightarrow F$. Then the orthogonal similitude group of $V$ with respect to the form $B$ is the group $\mathrm{GO}(V,B)$ of all $g \in \mathrm{GL}(V)$ such that $B(gv,gw)= \lambda (g)B(v,w)$ for any $v,w \in V$ with $\lambda (g) \in k^*$. The multiplicative character $\lambda : \mathrm{GO}(V,B) \rightarrow k^*$ is called the similitude character. Note that the orthogonal subgroup $\mathrm{O} (V,B)$ is equal to $Ker (\lambda ).$

    Suppose $V$ is a two dimensional vector space over $k$ with a symplectic form $\Sp $ defined by the determinant. That is to say $\Sp (v,w)=det(v,w)$ for any $v,w \in V$ which are expressed as column vectors with respect to a fixed base and $(v,w)$ is written as a $2 \times 2$ matrix. Then we can define a bilinear form $B$ on $V \otimes V$ by $B(v_1 \otimes w_1,v_2 \otimes w_2)=\Sp (v_1,v_2)\Sp (w_1,w_2)$. It is easy to check that $B$ is a non-degenerate and symmetric bilinear form on $V \otimes V$ and the image of the tensor product from $\mathrm{GL} (2,k) \times \mathrm{GL} (2,k)$ to $\mathrm{GL} (4,k)$ is a subgroup of $\mathrm{GO}(k^4,B_0)$ if we fix an isometry between $(V \otimes V, B)$ and $(k^4,B_0)$ where $B_0$ is the standard bilinear form of $k^4$ defined by $B_0(v,w)=v^tw$ for any $v,w \in k^4$. Therefore we have the following exact sequence,
    \[ 1 \rightarrow k^* \rightarrow \mathrm{GL}(2,k) \times \mathrm{GL}(2,k) \rightarrow \mathrm{GO}(k^4, B_0) \]
    since $\mathrm{GL}(2,k)$ is the symplectic similitude group of $V$ with similitude character $\lambda (g)=det(g)$. In particular, we have the following exact sequence:
    \[ 1 \rightarrow \{\pm(I_2,I_2)\} \rightarrow \mathrm{SL}(2,k) \times \mathrm{SL}(2,k) \rightarrow \mathrm{SO}(k^4, B_0) \]

    From the discussion on page 57 of \cite{D}, the abelianization of $\mathrm{SO}(k^4,B_0)$ is isomorphic to $k^*/k^{*2}$, which is trivial if $k= \bar{k}$. Therefore if we assume $k= \bar{k}$, then $\mathrm{SO}(k^4,B_0)$ is equal to its commutator subgroup. By the discussion on page 59 of \cite{D} $\mathrm{SO}(k^4,B_0)/ \{\pm \mathrm{I}_4\}$ is isomorphic to the group product $\mathrm{PSL}(2,k) \times \mathrm{PSL}(2,k)$. Thus the map from $\mathrm{SL}(2,k) \times \mathrm{SL}(2,k)$ to $\mathrm{SO}(k^4,B_0)$ is onto and the map from $\mathrm{GL}(2,k) \times \mathrm{GL}(2,k)$ to $\mathrm{GO}(k^4,B_0)$ is onto $\mathrm{GSO}(k^4,B_0)$.

    If we use this:
    \[ 1 \rightarrow \G _m \rightarrow \mathrm{GL}(2) \times \mathrm{GL}(2) \rightarrow \mathrm{GSO}(B_0) \rightarrow 1 \]
    and apply $H^i$, then we can get:
    \[ 1 \rightarrow k^* \rightarrow \mathrm{GL}(2,k) \times \mathrm{GL}(2,k) \rightarrow \mathrm{GSO}(k^4, B_0) \rightarrow 1 \]
    since $H^1(Gal(\bar{k}/k),~\bar{k}^*)=1$.

	\begin{lem}
	Let $k~$be a number field with algebraic closure $\bar{k}$. Then
	\[ \mathrm{GSO} (4,k)= \frac{\mathrm{GL} (2,k) \times \mathrm{GL} (2,k)} { \{(c\mathrm{I}_2,~c^{-1}\mathrm{I}_2) \}},\]
	where $c\in k^*$
	\end{lem}
	\textit{Proof}.
	Assume $k=\bar{k}$.

	We start with some notations. We let $\mathrm{B}(v,~w)={}^tvw~$be a non-degenerate symmetric bilinear form,
\[\mathrm{GO} (n,k)=\{g \in \mathrm{GL}(n,k)|\mathrm{B}
	(gv,~gw)=\lambda (g)\mathrm{B}(v,~w),~\lambda(g) \in k^*,~v,~w \in k^n\},\]
where the multiplicative character $\lz : \GO (n,k) \to k^*$ is called the similitude character,
\[\mathrm{O}(n,k)=\{g \in \mathrm{GL}(n,k)|\mathrm{B}(gv,gw)=\mathrm{B}(v,~w)\},\]
\[\mathrm{SO}(n,k)=\{g \in \mathrm{O}(n,k)|\mathrm{det}g=1\}\]
and
\[\mathrm{Z}(n,k)= \mathrm{center~of~} \mathrm{GO}(n,k).\]

	For all $g \in \mathrm{GO}(4,k)$, $(\mathrm{det}g)^2=\lambda(g)^4~$and
\[\mathrm{GSO}(4,k)=\{ g \in \mathrm{GO}(4,k) | \mathrm{det}g= \lambda (g)^2 \}\]
	is generated by $\mathrm{SO}(4,k)$, $\mathrm{Z}(4,k)~$and $\mathrm{SO}(4,k) \cap \mathrm{Z}(4,k)=\{ \pm \mathrm{I}_4 \}$.

	First, let $\mathrm{W}~$be $k^2~$with the standard symplectic form given by determinant. Then the induced bilinear form $\mathrm{B}_1~$on $\mathrm{W}
	\otimes \mathrm{W}~$is non-degenerated and symmetric, and $B_1~$is given by $B_1(v_1 \otimes w_1,v_2 \otimes w_2)=det(v_1,v_2)det(w_1,w_2)$. There is an
	isometry between $(\mathrm{W} \otimes \mathrm{W,~B_1})~$and $(k^4,~\mathrm{B})$.

	Since $\mathrm{GL}(2,k)~$is the symplectic similitude group of $(\mathrm{W},~\mathrm{det})$, we can get a sequence,
    \[ 1 \to  k^* \stackrel{\iota}{\rightarrow} \mathrm{GL}(2,k) \times \mathrm{GL}(2,k) \stackrel{\beta}{\rightarrow} \mathrm{GO}(4,k)\]
	in which the map $\iota~$is given by $\iota (c)=(c\mathrm{I}_2,~c^{-1}\mathrm{I}_2)$.

	Let $\beta:\mathrm{GL}(2,k) \times \mathrm{GL}(2,k) \to \mathrm{GO}(4,k)$ be defined as follows.
 The quadratic space
	$(k^4,~\mathrm{B})~$is isometric to $(\mathrm{M_2}(k),~\mathrm{B_2})~$where $\mathrm{B_2}~$is the symmetric bilinear map $(X,~Y) \to tr({}^tXY)$. Under
	this identification, $\beta(g_1,~g_2)~$is the automorphism of $k^4~$given by $X \to {}^tg_1Xg_2~$for all $g_1,~g_2 \in \mathrm{GL}(2,k)$. And ker
	$\beta =\{ (t\mathrm{I}_2,t^{-1}\mathrm{I}_2) | t \in k^*\}$.

	We can calculate det$\beta (g_1,~g_2)=[$det$(g_1)$det$(g_2)]^2,~ \lambda (\beta (g_1,~g_2))=$det$(g_1)$det$(g_2)$. So, det$\beta (g_1,~g_2)=
	\lambda (\beta (g_1,~g_2))^2$. Therefore, image of $\beta \subset \mathrm{GSO}(4)$.

	Since $\mathrm{GSO}(4,k)~$is generated by $\mathrm{SO}(4,k)~$and $\mathrm{Z}(4,k)~$and $\mathrm{SO}(4,k) \cap \mathrm{Z}(4,k)=\{ \pm I_4\}$, it is
	enough to show that $\mathrm{Z}(4,k)$ is contained in the image of $\beta~$and $\mathrm{SO}(4,k)$  is contained in the image of $\beta~$. First part is clear in case $k$ is an algebraically closed field and for the second part, we know
	\[\mathrm{SO}(4,k)= \beta (\mathrm{SL}(2,k) \times \mathrm{SL}(2,k)) \subset \beta (\mathrm{GL}(2,k) \times \mathrm{GL}(2,k)).\]

	So $\ds \frac{\mathrm{GL}(2,k) \times \mathrm{GL}(2,k)}{\mathrm{ker} \beta}=\frac{\mathrm{GL}(2,k) \times \mathrm{GL}(2,k)}{\{(t\mathrm{I}_2,~t^{-1}
	\mathrm{I}_2)\}}$ is the image of $\beta $ which is now $\mathrm{GSO}(4,k)$.

	We can also deduce the following exact sequences when $k$ is not an algebraically closed field:
	\[ 1 \to k^* \stackrel{\iota}{\rightarrow} \mathrm{GL}(2,k) \times \mathrm{GL}(2,k) \stackrel{\beta}{\rightarrow} \mathrm{GSO}(4,k) \to H^1(Gal(\bar{k}/k),~\bar{k}^*)=1.\]
	Therefore, $\ds \mathrm{GSO}(4,k)= \frac{\mathrm{GL}(2,k) \times \mathrm{GL}(2,k)}{\{(c\mathrm{I}_2,~c^{-1}\mathrm{I}_2)\}}$.$\square$

\bigskip

\item[{\large\bf{4}}  ] {\large\bf{Langlands Parameter of GSp(4)}}

\bigskip



    Let $k$ be a number field.




We have a commutative diagram:\newline

$\begin{array}[c]{ccccccccc}
0&{\rightarrow}&k^*&{\rightarrow}&\GL (2,k) \times \GL (2,k)&\stackrel{\beta}{\rightarrow}&\GSO (4,k)&{\rightarrow}&0\\
 &      &\uparrow&             &\uparrow              &                           &\uparrow&             & \\
0&{\rightarrow}&k^* &{\rightarrow}&k^* \times k^*&\stackrel{\alpha }{\rightarrow}&k^* &{\rightarrow}&0\\
&&&&&&&&
\end{array}$

where $\beta$ is defined in the previous chapter and $\alpha =\beta |_{ k^* \times k^*}$.


\begin{lem}\cite{HST}
	There is a bijection between cuspidal automorphic representations $\widetilde{\pi}$ of $\mathrm{GSO}(4,\mathbb{A})$ and pairs ($\pi,~ \widetilde{\chi}$)
	of a cuspidal automorphic representation $\pi$ of $\mathrm{GL}(2, \mathbb{A}) \times \mathrm{GL}(2, \mathbb{A})$ and a gr\"{o}ssencharacter $ \widetilde{\chi}: k^* \diagdown
	\mathbb{A}^* \to \mathbb{C}^* ~$such that $\widetilde{\chi} \circ \alpha$ is the central character of $\pi$.
	\end{lem}
	\textit{Proof}.
	Since the bijection sends
	$\widetilde{\pi}~$to $(\{f \circ \beta|f \in \widetilde{\pi}\},~\chi_{\widetilde{\pi}})$, where $\chi_{\widetilde{\pi}}~$denote the central character of
	$\widetilde{\pi}$ and $\beta$ is the natural map from $\mathrm{GL}(2) \times \mathrm{GL}(2)$ to $\GSO (4)$ as above. In the other direction it sends the pair ($\pi,~\widetilde{\chi}$) to the set of functions from $\mathrm{GSO}(4, k ) \diagdown \mathrm{GSO}(4, \mathbb{A} )~$to
	$\mathbb{C}~$such that $f \circ \beta \in \pi~$and $f(zg)= \widetilde{\chi}(z)f(g)~$for all $z \in \mathbb{A}^* ~$and all $g \in \mathrm{GSO}(4,\mathbb{A})$.$\square$
	
    Note that the second set in the lemma maps 2-1 to the set of cuspidal automorphic representations of $\mathrm{GL}(2, \mathbb{A}) \times \mathrm{GL}(2, \mathbb{A})$ whose central character factors through the map $\alpha$.
	Moreover, we can apply the same considerations for the local case. We consider the non-archimedean place $v$. There is a bijection between irreducible admissible representations $\widetilde{\pi}_v~$
	of $\mathrm{GSO}(4, k_v)~$and pairs ($\pi_v,~ \widetilde{\chi}_v$) of an irreducible representation $\pi_v~$of $\mathrm{GL}(2, k_v) \times \mathrm{GL}(2, k_v)$ and a character $ \widetilde{\chi}_v: k_v^* \to \mathbb{C}^*~$such that $\widetilde{\chi}_v \circ \alpha $ is the central character of $\pi_v$.


	

	
	Let $G$ be $\mathrm{GL}(2) \times \mathrm{GL}(2)$. For the rest of this section induction will mean unitary induction. Let $B_G~$denote the Borel subgroup of upper triangular matrices in $G$. Four
	characters $\chi_{11},\chi_{21},\chi_{12},\chi_{22}~$of $k_v^*~$give rise to a character $(\chi_{11},\chi_{21},\chi_{12},\chi_{22})~$of $B_G(k_v)~$by:
	
	\begin{displaymath}
	(\chi_{11},\chi_{21},\chi_{12},\chi_{22})
	\left( \begin{array}{cccc}
	 d_1 &  *  &     &     \\
	  0  & d_2 &     &     \\
         &     & d_3 & *   \\
         &     &  0  & d_4
	\end{array} \right)
	=\chi_{11}(d_1)\chi_{21}(d_2)\chi_{12}(d_3)\chi_{22}(d_4).
	\end{displaymath}

	We let $T_G~$denote the torus of diagonal matrices.

%

	Let $B_{\GO (4)}~$denote the Borel subgroup of $\GO (4)~$
    \begin{displaymath}
	B_{\GO(4)}=
	\left\{\left( \begin{array}{cccc}
	  a & * &   &    \\
	    & b &   &    \\
        &   & c & *  \\
        &   &   & d
	\end{array} \right)
    \in \GO(4)\right\}.
	\end{displaymath}

     Let $T_{\GO (4)}~$denote the Levi component
%


\begin{displaymath}
	T_{\GO(4)}=
	\left\{\left( \begin{array}{cccc}
	  a &   &   &   \\
	    & b &   &   \\
        &   & c &   \\
        &   &   & d
	\end{array} \right)
    \in \GO(4)\right\}.
	\end{displaymath}

%
	Let ($\pi, \wt{\chi}$) be a pair as Lemma 2 corresponding to $\wt{\pi}$. Suppose that $\pi_v~$is the principal series corresponding to a character
	$(\chi_{11},\chi_{21},\chi_{12},\chi_{22})~$of $B_G(k_v)$. Then $\chi_{11}\chi_{21}=\chi_{12}\chi_{22}=\wt{\chi}_v$, by page 384\cite{HST}.
	


Let $\mu~$and $\nu~$denote the multiplier characters of $\GSp(4)$ and $\GO (4)$ and let ${\mathrm{Sp(4)}}~$and
	$\mathrm{O}(4)~$(resp.) denote their kernels. Let $R=\{(g,h) \in \GSp (4) \times \GO (4) : \mu (g) \nu (h)=1\}$.
	
The group $Q$
is a minimal parabolic subgroup of $\GO (4)$.
Let $R_Q=R \cup (\GSp (4) \times Q)$.

The group
\begin{displaymath}
	P=
	\left\{\left( \begin{array}{cccc}
	  * & * & * & * \\
	    & * & * & * \\
        &   & * & *  \\
        &   &   & *
	\end{array} \right)
    \in \GSp(4)\right\}.
	\end{displaymath}
is a minimal parabolic subgroup of $\GSp (4)$.
Let $R_{P,Q}=R \cap (P \times Q)$.


	From a standard calculation, we can get following result on the Langlands parameters.
	\begin{lem}
	The $L-$group of $\GSp (4)~$is $\GSp (4, \C)$. If $\Pi~$is the unramified sub-quotient of the representation of $\GSp (4, k_v)~$unitarily induced from the
	character of $P(k_v)~$which is trivial on the unipotent radical and sends:
	\[diag(a,b,\mu a^{-1},\mu b^{-1}) \to \chi_1(a)\chi_2(b)\chi_3 (\mu),\]
	then $\Pi~$has Langlands parameter $(\chi_3(v),\chi_3\chi_1(v),\chi_2\chi_1(v),\chi_3\chi_1\chi_2(v)) \in \GSp (4, \C)$.
	\end{lem}	
	
	The following is from Rodier's classification which we need for the proof of proposition 2.
	\begin{lem}(Rodier's classification\cite{R})
	Suppose $\Pi~$is an irreducible pre-unitary representation of $\mathrm{GSp} (4, k _v)~$which is a subquotient of an unramified principal series
	representation with Langlands parameter $diag(\alpha,~\beta,~\gamma,~\delta) \in \mathrm{GSp} (4, \mathbb{C})$, then either $\Pi~$is the full induced
	representation or absolute value of $\alpha,~\beta,~\gamma,~\delta~$are, up to the action of the Weyl group, $\nu~$to the power
	$\ds (-\frac{1}{2},~-r,~r,~\frac{1}{2})~$with $\ds 0 \leq r \leq \frac{1}{4}$, or $\ds (-\frac{1}{2},~-\frac{1}{2},~\frac{1}{2},~\frac{1}{2})$,
	or ($\ds -\frac{3}{2},~-\frac{1}{2},~\frac{1}{2},~\frac{3}{2}$).
	\end{lem}

	The main proposition to get the Langlands parameter for $\GSp (4)~$ when the representation is associated to the representation of $\GL (2) \times
	\GL (2)$ is the following.
	\begin{prop}
	\cite{HST} Suppose that $\pi =\pi_1 \bigotimes \pi_2~$is an unramified irreducible pre-unitary principal series representation of
	$\mathrm{GL} (2) \times \mathrm{GL} (2)~$ with Langlands parameters diag($\alpha_1,~\beta_1$) and diag($\alpha_2,~\beta_2$).
	Suppose that $\Pi~$ is a pre-unitary irreducible admissible representation of $\mathrm{GSp} (4)~$which is associated to the representation ($\pi,~\tilde{\chi}$) obtained by theta lifting.
	Then $\Pi~$is an unramified irreducible principal series representation of $\mathrm{GSp} (4)~$with Langlands parameter diag($\alpha_1,~\alpha_2,~\beta_1,
	~\beta_2) \in \mathrm{GSp} (4, \mathbb{C}$).
	\end{prop}

	\textit{Proof}.
	The representation ${\pi}$ is induced from two pairs of characters ($\chi_{11},~\chi_{21}$) and ($\chi_{12},~\chi_{22}$) with
	$\tilde{\chi}=\chi_{11}\chi_{21}=\chi_{12}\chi_{22}$. Here, characters of the torus of $\mathrm{GSO}(4)$ are defined as\newline
	\begin{displaymath}
	{\chi_1(t_1,~t_2,~t_3):} =
	\left( \begin{array}{cccc}
	 t_1 &     &             &             \\
	     & t_2 &             &             \\
	     &     & t_3t_1^{-1} &             \\
	     &     &             & t_3t_2^{-1}
	\end{array} \right)
	{\mapsto (\frac{\chi_{11}}{\chi_{12}})(t_1)|t_2|(\frac{\chi_{12}}{\chi_{21}})(t_2)|t_3|^{-\frac{1}{2}}\chi_{21}(t_3)}
	\end{displaymath}

	or one of its conjugates under the group W of order 8 which is generated by $\sigma_1$, which switches $\chi_{11}~$and $\chi_{21}$, and $\tau~$which
	switches $\chi_{j1}~$and $\chi_{j2}~$for $j=1,~2$. Because $\pi~$is unitary and irreducible principle series, $\chi_{ij} \neq \chi_{i'j'}|~|$. Let
	$\mathrm{R}=ker\mu \nu$, where $\mu~$is the similitude character of $\mathrm{GSp} (4)~$and $\nu~$is the similitude character of
	$\mathrm{GO}(4)~$and $\mathrm{P}~$is the minimal parabolic subgroup of $\mathrm{GSp} (4)$, $\mathrm{Q}~$is the minimal parabolic subgroup of
	$\mathrm{GO}(4)~$.

	Therefore, for one of the characters $\chi_1~$above, $\Pi \bigotimes \chi_{1}~$must be a quotient of the induction from $\mathrm{R} \bigcap (\mathrm{P}
	 \times \mathrm{Q})~$to $\mathrm{R} \bigcap (\mathrm{GSp} (4) \times \mathrm{Q}$) of the character which is trivial on the unipotent radical and sends

    \begin{displaymath}
    \left( \begin{array}{cc}
	\left( \begin{array}{cccc}
	a& &          &          \\
	 &b&          &          \\
	 & &\mu a^{-1}&          \\
	 & &          &\mu b^{-1}
	\end{array} \right)
    ,
	\left( \begin{array}{cccc}
	t_1&   &              &              \\
	   &t_2&              &              \\
	   &   &(\mu t_1)^{-1}&              \\
	   &   &              &(\mu t_2)^{-1}
	\end{array} \right)
    \end{array} \right)
    \end{displaymath}

    \[\mapsto |\mu|^{-2}|ab|^2 \chi_1(a^{-1}t_1 \mu,~ b^{-1}t_2 \mu,~1)\]
	
	Therefore $\Pi \bigotimes \chi_1(1,~1,~\mu^{-1})~$ must be a quotient of the un-normalized induction from
 $\mathrm{P}~$ to $\mathrm{GSp} (4)~$ of a
	character which is trivial on unipotents and sends:

	\begin{displaymath}
	\left( \begin{array}{cccc}
	 a &   &            &            \\
	   & b &            &            \\
	   &   & \mu a^{-1} &            \\
	   &   &            & \mu b^{-1}
	\end{array} \right)
    {\mapsto |\mu^{-\frac{1}{2}}||b| \chi_1 (a^{-1}\mu,~b^{-1}\mu,~1)}
	\end{displaymath}

	for one of the characters $\chi_1$.

	Since \[\chi_1(1,~1,~\mu^{-1})=(\frac{\chi_{11}}{\chi_{12}})(1)|1|(\frac{\chi_{12}}{\chi_{21}})(1)|\mu^{-1}|^{-\frac{1}{2}}\chi_{21}(\mu^{-1})=
	|\mu|^{\frac{1}{2}}\chi_{21}(\mu^{-1}),\]
	\[ |\mu^{-\frac{1}{2}}||b| \chi_1 (a^{-1}\mu,~b^{-1}\mu,~1) \cdot |\mu|^{-\frac{1}{2}}\chi_{21}(\mu)\]
	\[ =|\mu|^{-\frac{1}{2}}|b|(\frac{\chi_{11}}{\chi_{12}})(a^{-1} \mu)|b^{-1} \mu|(\frac{\chi_{12}}{\chi_{21}})(b^{-1} \mu)|1|^{-\frac{1}{2}}
	\chi_{21}(1) \cdot  |\mu|^{-\frac{1}{2}}\chi_{21}(\mu)\]
	\[ =(\frac{\chi_{12}}{\chi_{11}})(a)(\frac{\chi_{21}}{\chi_{12}})(b) \chi_{11}(\mu).\]

	Therefore, $\Pi~$is a quotient of the un-normalized induction from $\mathrm{P}~$of the character which sends

	\begin{displaymath}
	\left( \begin{array}{cccc}
	 a &   &            &            \\
	   & b &            &            \\
	   &   & \mu a^{-1} &            \\
	   &   &            & \mu b^{-1}
	\end{array} \right)
	{\mapsto (\frac{\chi_{12}}{\chi_{11}})(a)(\frac{\chi_{21}}{\chi_{12}})(b)\chi_{11}(\mu)}
	\end{displaymath}

	or one of its conjugates by W.

	The un-normalized induction of all these characters have unramified subquotients with Langlands parameters

	\begin{displaymath}
	\left( \begin{array}{cccc}
	\chi_{11}(\nu) &                &                &                \\
	               & \chi_{12}(\nu) &                &                \\
	               &                & \chi_{21}(\nu) &                \\
	               &                &                & \chi_{22}(\nu)
	\end{array} \right)
	\in \mathrm{GSp} (4)
	\end{displaymath}
	
	If $\Pi~$is not the full induced representation, by Rodier's classification, $|\chi _{ij} (\nu)| = |\nu |^{\az _{ij}}~$with $\ds (\az _{11},\az _{12},
	\az _{21},\az _{22})=(-\frac{1}{2},~-r,~r,~\frac{1}{2})~$with $\ds 0 \leq r \leq \frac{1}{4}$, or
	$\ds (-\frac{1}{2},~-\frac{1}{2},~\frac{1}{2},~\frac{1}{2})$,
	or ($\ds -\frac{3}{2},~-\frac{1}{2},~\frac{1}{2},~\frac{3}{2}$). But since $|\chi _i (\nu)| < |\nu |^{\frac{1}{2}}$, $\Pi~$is full induced
	representation. Therefore, the result follows. $\square$

\bigskip

\item[{\large\bf{5}}  ] {\large\bf{Transfer from GSp(4) to GL(4)}}

\bigskip
	

This section is from \cite{AS}.

	Let $\mathbb{A}=\mathbb{A}_k~$denote the ring of adeles of a number field $k$. Let $\pi~$be a unitary cuspidal representation of $\mathrm{GSp}(4,\mathbb{A}_k)$, which
	we assume to be globally generic. Then $\pi~$has a unique transfer to an
	automorphic representation $\Pi~$of $\mathrm{GL}(4,\mathbb{A}_k)$. The transfer is
	generic (globally and locally) and satisfies $\omega _{\Pi}=\omega _{\pi}^2~$and
	$\Pi \simeq \widetilde{\Pi} \otimes \omega _{\pi}$. Here, $\omega _{\pi}~$and $\omega _{\Pi}~$denote the central characters of $\pi~$and $\Pi$, respectively.
	Moreover \cite{AS} gives a cuspidality criterion for $\Pi~$and proves, when $\Pi~$is not cuspidal, it is an isobaric sum of two unitary cuspidal representations of
	$\mathrm{GL}(2,\mathbb{A}_k)$.
	We define the similitude symplectic group of degree four via
	\[\GSp (4)=\{g \in \GL (4):{}^tgJg=\mu (g)J\},\]
	where
	\begin{displaymath}
	{J =}
	\left( \begin{array}{cccc}
	     &     &    &   1  \\
	     &     & 1  &      \\
	     &  -1 &    &      \\
	  -1 &     &    &
	\end{array} \right)
	\end{displaymath}
	and $\mu $ is the similitude character. We fix the following parametrization of the elements of the maximal torus $\bf {T}~$in
	$\mathrm{GSp}(4)$:
	\begin{displaymath}
	{\bf{T}} =
	\left\{t=t(a_0,a_1,a_2)=
	\left( \begin{array}{cccc}
	 a_0a_1a_2  &        &         &      \\
	            & a_0a_1 &         &      \\
	            &        & a_0a_2  &      \\
	            &        &         & a_0
	\end{array} \right)
	\right\}.
	\end{displaymath}

	Let $\pi= \otimes _v \pi _v~$be a globally $\psi$-generic unitary cuspidal automorphic representation of $\mathrm{GSp}(4,\mathbb{A})$. Here, $\psi= \otimes _v \psi
	_v~$is a non-trivial additive character of $k \setminus \mathbb{A}~$defining a character of the unipotent radical of the standard upper-triangular Borel subgroup in
	the usual way. We fix $\psi$ throughout this paper. Let $S~$be any non-empty finite set of non-archimedean places $v$, which
	includes those $v~$with $\pi _v~$or $\psi _v~$ramified. Asgari and Shahidi prove that there exists an automorphic representation $\Pi= \otimes _v \Pi _v~$of
	$\mathrm{GL}(4, \mathbb{A})~$such that $\Pi _v~$is a local transfer of $\pi _v~$for outside of $S$.

	To be more explicit, assume that $v \notin S$. If $v~$is archimedean, then $\pi _v~$is given by a parameter $\phi _v : W_v \to \mathrm{GSp}(4, \mathbb{C})$, where
	$W_v~$is the Weil group of $k_v$. Let $\Phi _v : W_v \to GL(4,\mathbb{C})~$be given by $\Phi _v = \iota \circ \phi _v$, where $ \iota : \mathrm{GSp}(4,\mathbb{C}) \to
	\mathrm{GL}(4,\mathbb{C})~$is the natural embedding. Then $\Phi _v~$is the parameter of $\Pi _v$.
	
	If $v \notin S~$is non-archimedean, then $\pi_v~$is the unique unramified subquotient of the representation induced from an unramified character $\chi~$of $\bf
	{T}$$(k_v)~$to $\mathrm{GSp}(4,k_v)$. Writing $\chi(t(a_0,a_1,a_2))= \chi _0(a_0)\chi _1(a_1)\chi _2(a_2)$, where $\chi _i~$are unramified characters of $k_v^{\times}~$
	and $a_i \in k_v^{\times}$, the representation $\Pi _v~$is then the unique irreducible unramified subquotient of the representation of
	$\mathrm{GL}(4,k_v)~$parabolically induced from the character
	\[\chi_1 \otimes \chi_2 \otimes \chi_2^{-1} \chi_0 \otimes \chi_1^{-1}\chi_0~\]
	of $\bf{T}$$(k_v)$.
	
	Moreover, they proved that $\oz _ {\Pi} =\oz ^2$, where $\oz =\oz _ {\pi}~$and $\oz _{\Pi}~$denote the central characters of $\pi~$and $\Pi$, respectively, and
	for $v \not\in S~$they have $\Pi _v \sim \wt{\Pi} _v \otimes \oz _{\pi_v}$, i.e. $\Pi~$is nearly equivalent to $\wt{\Pi} \otimes \oz$.
	
	The representation $\Pi~$equivalent to a subquotient of some representation
	\[Ind(|\mathrm{det}|^{r_1}\sigma _1 \otimes \dots \otimes |\mathrm{det}|^{r_t}\sigma _t),\]
	where induction is from $\mathrm{GL}(n_1) \times \cdots \times
	\mathrm{GL}(n_t)~$with $n_1+ \cdots +n_t=4~$to $\mathrm{GL}(4)~$and $\sigma _i~$are the unitary cuspidal automorphic representation of
	$\mathrm{GL}(n_i,~\mathbb{A})~$and $r_i \in \mathbb{R}$.
	
	Without loss of generality we may assume that $r_1 \geq r_2 \geq \cdots \geq r_t$. Moreover, as $\Pi~$is unitary we have $n_1r_1+ \cdots +n_tr_t=0$, which
	implies that $r_t \leq0$. Let $T=S \cup \{v:v| \infty\}~$and consider
	\[L^T(s, \pi \times \wt{\sigma}_t)=L^T(s, \Pi \times \wt{\sigma}_t)= \prod _{i=1}^t L^T(s+r_i,\sigma_i \times \wt{\sigma}_t).\]
	Here, $L^T~$denotes the product over $v \not\in T~$of the local L-functions.
	
	If $n_t=1$, then the left-hand side is entire by a result of Piatetski-Shapiro \cite{P}. Now consider the right-hand side at $s_0=1-r_t \geq 1$. The last term in the
	product has a pole at $s_0$, whereas all of the others are non-zero there as $R(s_0 + r_i)=1+r_i-r_t \geq 1$. This is a contradiction.
	
	Now assume that $n_t=3$, i.e. $t=2~$with $n_1=1~$and $n_2=3$. Replacing $\pi~$and $\Pi$ by their contragredients will change $r_i~$to $-r_i~$and takes us back
	to the above situation, which gives a contradiction again.
	
	Therefore, $n_t=2$. In this case, $L^T(s, \pi \times \wt{\sigma}_t)~$have a pole at $s=1~$and if so, arguing as above, we conclude that $r_t=0$. This means
	that we either have $t=2~$with $n_1=n_2=2~$or $t=3~$with $n_1=n_2=1~$and $n_3=2$. However, we can rule out the latter as follows.
	
	Assume that $t=3~$with $n_1=n_2=1~$and $n_3=2$. Then, it follows from the fact that $r_3=0~$and contradictions $r_1 \geq r_2 \geq r_3~$and $r_1+r_2+2r_3=0~$
	that all of the $r_i~$would be zero in this case. This implies that if we consider the L-function of $\pi~$twisted by $\wt{\sigma}_1$, we have
	\[L^T(s, \pi \times \wt{\sigma}_1)=L^T(s, \sigma _1 \times \wt{\sigma}_1) L^T(s, \sigma _2 \times \wt{\sigma}_1) L^T(s, \sigma _3 \times \wt{\sigma}_1).\]
	Now the left-hand side is again entire by Piatetski-Shapiro's result \cite{P} mentioned above and the right-hand side has a pole at $s=1$, which is a contradiction.
	
	Therefore, the only possibilities are $t=1$ i.e. $\Pi~$unitary cuspidal or $t=2~$and $n_1=n_2=2~$with $r_2=0$. In the latter case, we also get $r_1=0$, as
	$r_1+r_2=0~$by unitarity of the central character. Moreover, in this case we have $\sigma _1 \not\simeq \sigma _2~$as, otherwise,
	\[L^T(s,\pi \times \wt{\sigma _1})=L^T(s, \sigma_1 \times \wt{\sigma_1}) L^T(s, \sigma_2 \times \wt{\sigma_1})\]
	must have a double pole at $s=1~$while any possible pole of the left-hand side at $s=1~$is simple.
	
	Therefore, we can see the following.
	
	\begin{prop}
	
	\cite{AS} Let $\pi~$be globally generic unitary cuspidal automorphic representation of $\mathrm{GSp}(4, \mathbb{A})~$and let $\Pi~$be any transfer of $\pi~$to
	$\mathrm{GL}(4,\mathbb{A})$. Then $\Pi~$is a subquotient of an automorphic representation as Ind=($|\mathrm{det}|^{r_1}\sigma _1 \otimes \dots \otimes
	|\mathrm{det}|^{r_t}\sigma _t$) with either $t=1,~n_1=4~$and $r_1=0~$(i.e. $\Pi~$is unitary cuspidal) or $t=2,~n_1=n_2=2~$and $r_1=r_2=0$. In the latter case,
	 we have $\sigma _1 \not\simeq \sigma _2$.
	
	\end{prop}
	
	In fact, we can get more precise information.
	
	\begin{prop}
	
	\cite{AS} Let $\pi~$be a globally generic unitary cuspidal automorphic representation of $\mathrm{GSp}(4,\mathbb{A})~$with $\omega=\omega _{\pi}~$its central character and let
	$\Pi~$be any transfer as above. Then, $\Pi \simeq \wt{\Pi} \otimes \omega~$(not just nearly equivalent). Moreover:\newline
	(a) the representation $\Pi~$is cuspidal if and only if $\pi~$is not obtained as a Weil lifting from $\mathrm{GSO}(4,\mathbb{A})$\newline
	(b) if $\Pi~$is not cuspidal, then it is the isobaric sum of two representations $\Pi= \Pi _1 \boxplus \Pi_2$, where each $\Pi _i~$is a unitary cuspidal
	automorphic representation of $\mathrm{GL}(2,\mathbb{A})~$satisfying $\Pi _i \simeq \widetilde{\Pi _i} \otimes \omega~$and $\Pi _1 \not\simeq \Pi_2$.
	
	\end{prop}

\bigskip

\item[{\large\bf{6}}  ] {\large\bf{L-functions}}

\bigskip


	By the natural embedding from $\mathrm{GSp} (4, \mathbb{C})~$to $\mathrm{GL} (4, \mathbb{C})$, we can see a representation $\Pi~$of $\mathrm{GL} (4)$ which
	is transferred from $\mathrm{GSp} (4)~$is not cuspidal when it is obtained as a Weil lifting from $\mathrm{GSO} (4, \mathbb{A} )$, and in this case it is
	the isobaric sum of two representations $\Pi_i$'s, where each $\Pi_i~$is a unitary cuspidal automorphic representation of $\mathrm{GL} (2, \mathbb{A})~$
	satisfying $\Pi_i\simeq \widetilde{\Pi_i}\otimes \omega~$and $\Pi_1 \not\simeq \Pi_2$.

%
	   \begin{thm}
	Let $\pi _i,~i=1,~2$, be cuspidal generic representations of $\mathrm{GSp}(4,\A)$ and $\Pi_i,~i=1,~2$, be their transfers.
	
	\begin{itemize}
	\item[1.]
	If neither of $\pi _i,~i=1,~2$ come from $\GSO (4, \A )$, then $L^S(s,\pi_1 \times \pi_2)$ has a pole at
	$s=1$ if and only if $\pi_2 = \widetilde{\pi _1}$.
	\item[2.]
	If only one of $\pi_i,~i=1,~2$ comes from $\mathrm{GSO}(4,\mathbb{A})$, then $L^S(s,\pi_1 \times \pi_2)$ has no poles.

	\item[3.]
	Suppose the representations $\pi_1$, $\pi_2~$of $\mathrm{GSp}(4,\A)$ are obtained as a Weil lifting from $\mathrm{GSO}(4,\mathbb{A} )$.
	Then
	$\Pi_1= \Pi_{11} \boxplus \Pi_{12}~$and $\Pi_2= \Pi_{21} \boxplus \Pi_{22}~$, and
	\[L^S(\Pi_1 \times \Pi_2)=L^S((\Pi_{11} \boxplus \Pi_{12}) \times (\Pi_{21} \boxplus \Pi_{22}))\]
	\[=L^S(\Pi_{11} \times \Pi_{21})L^S(\Pi_{11} \times \Pi_{22})L^S(\Pi_{12} \times \Pi_{21})L^S(\Pi_{12} \times \Pi_{22})\]

	Consequently,
	\begin{itemize}
	\item[(a)] if $\Pi_{11} \not\simeq \widetilde{\Pi}_{21}~$and $\Pi_{12}
	\not\simeq \widetilde{\Pi}_{22}$, then $L_s(\Pi_1 \times \Pi_2)~$has no poles, since $\Pi_{11} \not\simeq \Pi_{12},~\Pi_{21} \not\simeq \Pi_{22}$.
	\item[(b)] if $\Pi_{11} \simeq \widetilde{\Pi}_{21}~$and $\Pi_{12} \not\simeq \widetilde{\Pi}_{22}~$, then $L^S(\Pi_1 \times \Pi_2)~$has a simple pole at
	$s=1$.
	\item[(c)] if $\Pi_{11} \simeq \widetilde{\Pi}_{21}~$and $\Pi_{12} \simeq \widetilde{\Pi}_{22}~$, then $L^S(\Pi_1 \times \Pi_2)~$has a double
	pole at $s=1$.
	\end{itemize}
    \end{itemize}
	\end{thm}
	
	\textit{Proof.}
	We know that $L^S(\Pi \times \wt {\Pi})~$has a simple pole at $s=1~$when $\Pi~$is the representation of $\GL (2, \A)~$and from the last section we know
	that $\Pi_{11} \not\simeq \Pi_{12},~\Pi_{21} \not\simeq \Pi_{22}$.
	
	If the representation $\Pi_1~$is not obtained as a Weil lifting from $\mathrm{GSO} (4, \mathbb{A} )~$and $\Pi_2~$is obtained as a Weil lifting from
	$\mathrm{GSO} (4, \mathbb{A})$, then $\Pi_2= \Pi_{21} \boxplus \Pi_{22}$. And
	\[L^S(\Pi_1 \times \Pi_2)=L^S(\Pi_1 \times (\Pi_{21} \boxplus \Pi_{22}))=
	L^S(\Pi_1 \times \Pi_{21})L^S(\Pi_1 \times \Pi_{22})\]
	and we can see this L-function has no poles (c.f. \cite{JS}).

	If representations $\Pi_1$, $\Pi_2~$of $\mathrm{GL} (4,\A)$ are obtained as Weil liftings from $\mathrm{GSO} (4, \mathbb{A} )~$, then
	$\Pi_1= \Pi_{11} \boxplus \Pi_{12}~$and $\Pi_2= \Pi_{21} \boxplus \Pi_{22}~$, and
	\[L^S(\Pi_1 \times \Pi_2)=L^S((\Pi_{11} \boxplus \Pi_{12}) \times (\Pi_{21} \boxplus \Pi_{22}))\]
	\[=L^S(\Pi_{11} \times \Pi_{21})L^S(\Pi_{11} \times \Pi_{22})L^S(\Pi_{12} \times \Pi_{21})L^S(\Pi_{12} \times \Pi_{22})\]
	
	Therefore, if $\Pi_{11} \not\simeq \widetilde{\Pi}_{21}~$and $\Pi_{12} \not\simeq \widetilde{\Pi}_{22}$, then $L^S(\Pi_1 \times \Pi_2)~$has no poles,
	since $\Pi_{11} \not\simeq \Pi_{12},~\Pi_{21} \not\simeq \Pi_{22}$. If $\Pi_{11} \simeq \widetilde{\Pi}_{21}~$and $\Pi_{12} \not\simeq
	\widetilde{\Pi}_{22}~$, then $L^S(\Pi_1 \times \Pi_2)~$has a simple pole at $s=1$ because $L^S(\Pi_{11} \times \Pi_{21})~$has a simple pole and
	$L^S(\Pi_{11} \times \Pi_{22})$, $L^S(\Pi_{12} \times \Pi_{21})~$and $L^S(\Pi_{12}\times \Pi_{22})~$have no poles.
	
	If $\Pi_{11} \simeq \widetilde{\Pi}_{21}~$and $\Pi_{12} \simeq \widetilde{\Pi}_{22}~$, then $L^S(\Pi_1 \times \Pi_2)~$has a double
	pole at $s=1~$because $L^S(\Pi_{11} \times \Pi_{21})$, $L^S(\Pi_{12}\times \Pi_{22})$ each have a simple pole and
	$L^S(\Pi_{11} \times \Pi_{22})$, $L^S(\Pi_{12} \times \Pi_{21})~$have no poles and are non-zero at $s=1$ \cite{S81}.

Since $\Pi _1 \not\simeq \Pi_2$ by Proposition 4, above cases are all for this theorem.
	
	In fact, $L^S(\Pi \times \wt {\Pi})~$has a double pole at s=1 if $\Pi~$is a Weil lifting from $\GSO (4)$. $\square$		

    Thus part (c) shows the existence of a L-function of a cuspidal representation of $ \GSp (4,\A) \times \GSp (4,\A)~$which has a pole of order
	2 at
	$s=1$, even for globally generic representations.
	
	\begin{thm}
If $\pi$ comes from $\GSO (4,\A )$, then $\pi$ is the Weil transfer of $\Pi _1 \otimes \Pi _2$ realized as a representation of $\GSO (4,\A)$. This agrees with Langlands Functoriality principle as $\GSO (4)$ is an endoscopic group for $\GSp (4)$. Moreover, it shows the data $\Pi _1 \otimes \Pi_2$ on $\GSO (4)$ transfers to $\Pi_1 \boxplus \Pi_2$ through the composite of the endoscopic transfer from $\GSO (4)$ to $\GSp (4)$ and the twisted endoscopic transfer from $\GSp (4)$ to $\GL (4)$.
\end{thm}

\textit{Proof.} By Lemma 2, we know that there is a bijection between cuspidal automorphic representations $\widetilde{\Pi}$ of $\mathrm{GSO}(4,\mathbb{A})$ and pairs ($\Pi,~ \widetilde{\chi}$) of a cuspidal automorphic representation $\Pi$ of $\mathrm{GL}(2, \mathbb{A}) \times \mathrm{GL}(2, \mathbb{A})$ and a gr\"{o}ssencharacter $ \widetilde{\chi}: k^* \diagdown \mathbb{A}^* \to \mathbb{C}^* ~$such that $\widetilde{\chi} \circ \alpha$ is the central character of $\Pi$. By Lemma 1, we see $\GSO (4)$ is basically $\GL (2)\times \GL (2)$. Let $\Pi =\Pi_1 \bigotimes \Pi_2~$be an unramified irreducible pre-unitary principal series representation of $\mathrm{GL} (2) \times \mathrm{GL} (2)~$ with Langlands parameters diag($\alpha_1,~\beta_1$) and diag($\alpha_2,~\beta_2$). Then by Proposition 2, if we say $\pi~$ is a pre-unitary irreducible admissible representation of $\mathrm{GSp} (4)~$which is associated to the representation ($\Pi,~\tilde{\chi}$) obtained by the theta lifting, then $\pi~$is an unramified irreducible principal series representation of $\mathrm{GSp} (4)~$with Langlands parameter diag($\alpha_1,~\alpha_2,~\beta_1,~\beta_2) \in \mathrm{GSp} (4, \mathbb{C}$).

Since $\pi $ is obtained as a Weil lifting from $\GSO (4, \A )$, by Proposition 4, we know any transfer $\Pi '$ from $\pi $ is cuspidal and the isobaric sum of two representations $\Pi '=\Pi _1' \boxplus \Pi _2'$, where each $\Pi _i'$ is a unitary cuspidal automorphic representation of $\GL (2, \A )$.

From the classification theorem, we can say $\Pi _i=\Pi _i'$ for $i=1,~2$ after reordering if it is necessary. Therefore $\Pi =\Pi '=\Pi_1 \boxplus \Pi_2$ and the result follows. $\square$

\end{enumerate}


\begin{thebibliography}{99}



    \bibitem{AS} Asgari, M. and Shahidi, F.: Generic Transfer from $\GSp(4)~$to $\GL(4)$, Composito Math. 142(2006) 541-550

    \bibitem{B76} Borel, A.: Admissible representations of a semi-simple group over a local field with vectors fixed under an Iwahori subgroup, Inv. Math. 35 (1976), pp. 233-259.


    \bibitem{D} Jean A Dieudonn\'{e}. La g\'{e}om\'{e}trie des groupes classiques, Springer-Verlag, Berlin, 1971. Troisi\'{e}me \'{e}dition, Ergebnisse des Mathematik und ihrer Grenzgebiete, Band 5.





    \bibitem{HST} Harris, M. Soudry, D. Taylor, R.: l-adic representations associated to modular forms over imaginary quadratic fields, Invent. math. 112,
    377-411(1993)

    \bibitem{J} Jiang, D.:The degree 16 standard L-function of $\GSp (2)\times \GSp (2)$. Mem. Amer. Math. Soc. 588 (1996)


    \bibitem{JS} Jacquet,H. and Shalika,J. On Euler Products and the Classification of Automorphic Forms, Amer.J.Math 103(1981) 777-875




    \bibitem{P} Piatetski-Shapiro, I. I.: L-functions for $\GSp _4$. Olga Taussky-Todd: in memoriam. Pacific J. Math. (1997), Special Issue, 259-275.

    \bibitem{R} Rodier, F.: Sur les Representations non Ramifiees des Groupes Reductifs p-adiques; l'Exemple de GSp(4), Bull. Soc. Math. Fr. 116, 15-42 (1988)

    \bibitem{S81} Shahidi, F.: On certain L-functions. Amer. J. Math., 103(2):297-355, 1981.



\end{thebibliography}
\end{document}